\documentclass[12pt, a4paper]{article}
\usepackage{amsmath, amssymb, latexsym}
\renewcommand{\baselinestretch}{1.15}
 \newtheorem{theorem}{Theorem}[section]
\newtheorem{lemma}[theorem]{Lemma}

\newtheorem{corollary}[theorem]{Corollary}
\newtheorem{remark}[theorem]{Remark}

\renewcommand{\title}[1]
{\thispagestyle{empty}
\begin{center}
{\Large \bf #1}
\end{center}}

\newcommand{\authors}[1]
{\begin{center}
\renewcommand{\thefootnote}{\fnsymbol{footnote}}
\setcounter{footnote}{3} {\sc #1 }
\end{center}}

\newcommand{\ack}[1]{\footnote{#1}}

\newcommand{\address}[1]
{\vskip 5ex
\renewcommand{\baselinestretch}{1}
\footnotesize \normalsize
 #1 \\
}

\begin{document}

\title{Hardy's inequalities for monotone functions on partially ordered measure spaces}

\authors{
Nicola Arcozzi\ack{Research partially supported by
the COFIN project Harmonic Analysis of the Italian Minister
of Research.},  
Sorina Barza, 
Josep L. Garcia-Domingo\ack{Research partially supported by Grants MTM2004-02299 and  
2001SGR00069.}, 
and
Javier Soria\ack{Research partially supported by Grants MTM2004-02299 and  
2001SGR00069.\\{\sl  Keywords:} Hardy operator, weighted inequalities.\\{\sl  MSC2000:} 46E30, 46B25.} }

\bigskip

 {\narrower\noindent \textbf{Abstract.} \small{We characterize the weighted Hardy's inequalities for monotone functions in ${\mathbb R^n_+}.$  In dimension $n=1$, this recovers the classical theory of $B_p$ weights. For $n>1$, the result was only known for the case $p=1$. In fact, our main theorem is proved in the more general setting of partially ordered measure spaces.}\par}

\bigskip
 
\section{Introduction}
The theory of weighted inequalities for the Hardy operator, acting on monotone functions in ${\mathbb R_+}$, was first introduced in \cite{AM}. Extensions of these results to higher dimension have been considered only  in very specific cases. In particular, in the diagonal case, only for $p=1$ (see \cite{BaPeSo1}). The main difficulty in this context is that the level sets of the monotone functions are not totally ordered, contrary to the one-dimensional case where one considers intervals of the form $(0,a)$, $a>0$. It is also  worth to point out that, with no monotonicity restriction, the boundedness of the Hardy operator is only known in dimension $n=2$ (see \cite{Saw}, \cite{KP}, and also \cite{BHP} for an extension in the case of product weights).

In this work we completely characterize  the weighted Hardy's inequalities for all values of $p>0$, namely, the boundedness of the operator:

$$
S: L^p_{\rm dec}(u)\longrightarrow L^p(u),
$$ 

\noindent
where

$$
Sf(s,t)=\frac1{st}\int_0^s\int_0^tf(x,y)\,dydx,
$$

\noindent
and $L^p_{\rm dec}(u)$ is the cone of positive and decreasing functions, on each variable, in $L^p(u)=L^p(\mathbb{R}^2_+,u(x)\,dx)$ (we consider, for  simplicity, $n=2$, although the result holds in any dimension).

The techniques we are going to use were introduced in \cite{CaSo}, for the one-dimensional case, and apply also to a more general setting, which we now define:

\medskip

We will consider a  family of $\sigma$-finite measure spaces $(X,\mu_x)$ (where $\mu_x$ is a measure on $X$, for each $x\in X$), with a partial order $\le $ satisfying:

\medskip

\begin{enumerate}

\item If $X_x:=\{u\in X:u\le x\}$, then $\le$ restricted to $X_x$ is a total order. 

\item If $D$ is a decreasing set, with respect to the order $\le$ (i.e.,  $\chi_D$ is a decreasing function),  then $D$ is measurable.

\item  $\mu_x(X_x)=1$ (observe that $X_x$ is a decreasing set).

\item If $u\in X_x$, then $d\mu_x(y)=\mu_x(X_u) d\mu_u(y)$. In particular, $\mu_x(X_u)\mu_u(X_x)=1$.

\end{enumerate}

\bigbreak
\noindent The main examples are: 

\begin{itemize}
\item{}   $X={\mathbb R}_+$ with the usual order, and  $\mu_x(E)=x^{-1}|E|$.  This is the case considered in \cite{AM}.

\item{} $X$ is a tree with the usual order on geodesics, and $\mu_x(E)={\rm Card}(E)/|x|$, where $|x|={\rm Card}([o,x])$ and $[o,x]$ is the geodesic path joining the origin $o$ with any vertex $x$ of the tree. For more information on this case, see \cite{GS} and the references quoted therein. 

\item{}   ${\mathbb R}^2_+$ with the order given by $(a_1,b_1)\le(a_2,b_2)$ if and only if, $a_1=a_2$ and $b_1\le b_2$ (we could also choose to fix the second coordinate). For $x=(a,b)\in{\mathbb R}^2_+$, 

$$
\mu_x(E)=b^{-1}\int_{E\cap(\{a\}\times{\mathbb R}_+)}\,dt.
$$

\item{}  In many cases, we can easily get the existence of the family of measures $\mu_x$ by taking a non-negative measure $\mu$ on $X$, and defining $\mu_x= {\mu}/{\mu(X_x)}$.
\end{itemize}

\medbreak

We now define the Hardy operator as follows:
$$
Sf(x)=\int_{X_x} f(u)\,d\mu_x(u).
$$

This definition is similar to the one considered in \cite{BaPeSo1}. For the case of ${\mathbb R}_+$, 
$$
Sf(x)=\frac1x\int_0^x f(t)\,dt.
$$

\noindent
 On a tree, 
 $$
 Sf(x)=\sum_{y\in[o,x]} \frac{f(y)}{|x|},
 $$
 and for ${\mathbb R}^2_+$, 
 $$
 Sf(a,b)=\frac1b\int_0^b f(a,t)\,dt.
 $$

\medskip

One of the main techniques we are going to use is the following lemma. This is a kind of integration by parts:

\medskip
\begin{lemma}
\label{IP} Let $(X,\mu,\le)$ be a finite measure space with a total order $\le$, and $\alpha\in {\mathbb R}$. Then, there exists a constant $C_{\alpha}$ which depends only on $\alpha$ (and not on $(X,\mu,\le)$), such that
$$\bigg(\int_Xd\mu\bigg)^{\alpha}\le C_{\alpha}\int_X\bigg(\int_{X_u}d\mu\bigg)^{\alpha-1}\,d\mu(u).
$$
\end{lemma}
\medskip

\noindent{\bf Proof:} 
Since $\mu(X)<\infty$, dividing both sides by $(\mu(X))^{\alpha}$, it suffices to show that 
$$1\le C_{\alpha}\int_X\varphi^{\alpha-1}(u)\, d\mu(u),
$$
 where $\varphi(u)=\mu(X_u)$, and $\mu(X)=1$.

If $\alpha\le 1$, using that $0\le \varphi\le 1$, then $\varphi^{\alpha-1}(u)\ge 1$, and
$$
\int_X\varphi^{\alpha-1}(u)\, d\mu(u)\ge \mu(X)=1.
$$

If $1<\alpha\le 2$, then $\varphi^{\alpha-1}(u)\ge \varphi(u)$, and hence it suffices to prove it for $\alpha=2$.

If $2<\alpha$, using Jensen's inequality:
$$
\bigg(\int_X\varphi(u)\,d\mu(u)\bigg)^{\alpha-1}\le\int_X\varphi^{\alpha-1}(u)\,d\mu(u),
$$
and as before, it reduces to the case $\alpha=2$.

Finally, if $\alpha=2$, 
\begin{align*}
\int_X\bigg(\int_{X_u}\,d\mu(x)\bigg)\,d\mu(u)&=\int_X\bigg(\int_{\{x\le u\}}\,d\mu(u)\bigg)\,d\mu(x)\\
&=\int_X\bigg(1-\int_{\{x>u\}}\,d\mu(u)\bigg)\,d\mu(x)\\
&=1-\int_X\int_{\{u\le x\}}\,d\mu(u)\,d\mu(x)+\int_X\int_{\{u=x\}}\,d\mu(u)\,d\mu(x),
\end{align*}
(here we used that the order is total). Thus,
$$
\int_X\bigg(\int_{X_u}\,d\mu(x)\bigg)\,d\mu(u)\ge \frac12.
$$
$\hfill\Box$\bigbreak

\section{Weighted Hardy's inequality}

In this section we will prove the main theorem. In order to include all the examples, we need to consider a second weaker order $\prec$ satisfying:

\begin{enumerate}

\item If $x\le y$, then $x\prec y$.

\item If $f$ is $\prec$-decreasing, then $Sf$ is $\prec$-decreasing.

\end{enumerate}

 We recall that we still keep $\le$ to define the operator $S$.
 
 \medskip

\begin{remark}{\rm\ 

\begin{itemize}
\item{}  We can (and will in some cases) take $\prec$ to be $\le$. In fact, we only need to check that the second condition holds for $\le$: If $f$ is $\le$-decreasing, and $u\le x$, then $Sf(x)\le Sf(u)$ if and only if, 

$$
\int_{X_u}f(y)[1-\mu_x(X_u)]\,d\mu_u(y)-\int_{X_x\setminus X_u}f(y)\mu_x(X_u)\,d\mu_u(y)\ge 0,
$$ 

\noindent (here we have used that $d\mu_x(y)=\mu_x(X_u)d\mu_u(y)$), and this follows from the fact that $\inf f|_{X_u}\ge\sup f|_{X_x\setminus X_u}$.

\item If a function $f$  or a set $D$ are $\prec$-decreasing they are also $\le$-decreasing.

\item The main example we have in mind is  $\le$ in ${\mathbb R}^2_+$ as before, and $\prec$ the order given by the rectangles (which is clearly a weaker order): $(a_1,b_1)\prec (a_2,b_2),$ if and only if $a_1\le a_2$ and $b_1\le b_2$ (i.e., the rectangle in $\mathbb{R}^2_+$ determined by the origin and $(a_1,b_1)$ is contained in the one determined by the origin and $(a_2,b_2)$). To show the second condition, assume that $f$ is a function decreasing on each variable. Then, it is obvious that $y^{-1}\int_0^y f(x,t)\, dt$ is also decreasing on each variable.

\item We denote by $L^p_{\prec}(d\nu)$ the class of $\prec$-decreasing functions in $L^p(d\nu)$. As a general assumption, we will only consider cases for which measurability of the functions involved always holds.

\end{itemize}}
\end{remark}

\medskip

\begin{theorem} \label{genth} Let $(X,\mu_x,\le)$ and $\prec$ satisfy the conditions given above, for all $x\in X$. Let $d\nu$ be a measure on $X$, and $p>0$. Then, the Hardy operator is bounded
$$S: L^p_{\prec}(d\nu)\rightarrow L^p(d\nu),
$$
if and only if, there exist a constant $C>0$ such that, for all $\prec$-decreasing sets $D$,
\begin{align}\label{gencon}
\int_{X\setminus D}\mu_x^p(D\cap X_x)\,d\nu(x)\le C\nu (D).
\end{align}
\end{theorem}
\medskip

\noindent{\bf Proof:}  Consider $f=\chi_D$, where $D$ is $\prec$-decreasing. Then $(Sf)^p(x)=\mu_x^p(D\cap X_x)$, and hence
\begin{align*}
\Vert Sf\Vert_{L^p(d\nu)}^p&=\int_X\mu_x^p(D\cap X_x)\,d\nu(x)\\
&=\int_D\mu_x^p(D\cap X_x)\,d\nu(x)+\int_{X\setminus D}\mu_x^p(D\cap X_x)\,d\nu(x)\\
&\le C\nu(D).
\end{align*}
Thus, 

$$\int_{X\setminus D}\mu_x^p(D\cap X_x)\,d\nu(x)\le C\nu(D).
$$ 

Observe that since $D$ is also $\le$-decreasing, if $x\in D$, then $X_x\subset D$, and so $\mu_x^p(D\cap X_x)=\mu_x^p( X_x)=1$. Therefore, $\int_D\mu_x^p(D\cap X_x)\,d\nu(x)=\nu(D)$.

\medskip

Conversely, if $p\ge 1$ and $f\in L^p_{\prec}(d\nu)$, using the lemma with $(X_x,\mu,\le)$ and $d\mu(u)=f(u)d\mu_x(u)$, we have that, for a constant $C$ which does not depend on either $f$ or $x$,
\begin{align*}
(Sf)^p(x)&=\bigg(\int_{X_x}f(u)\,d\mu_x(u)\bigg)^p\le C \int_{X_x}\bigg(\int_{X_u}f(y)\,d\mu_x(y)\bigg)^{p-1}f(u)\,d\mu_x(u)\\
&=C \int_{X_x}\bigg(\int_{X_u}f(y)\,d\mu_u(y)\bigg)^{p-1}f(u)\mu_x^{p-1}(X_u)\,d\mu_x(u)\\
&=C \int_0^{\infty}\int_{\{g>t\}\cap X_x}\mu_x^{p-1}(X_u)\,d\mu_x(u)\,dt,
\end{align*}
where $g(u)=(Sf)^{p-1}(u)f(u).$ Hence, 
\begin{align*}
\Vert Sf\Vert_{L^p(d\nu)}^p&\le C\int_X\int_0^{\infty}\int_{\{g>t\}\cap X_x}\mu_x^{p-1}(X_u)\,d\mu_x(u)\,dt\,d\nu(x)\\
&\approx\int_X\int_0^{g(x)}\int_{\{g>t\}\cap X_x}\mu_x^{p-1}(X_u)\,d\mu_x(u)\,dt\,d\nu(x)\\
&\qquad+
\int_X\int_{g(x)}^{\infty}\int_{\{g>t\}\cap X_x}\mu_x^{p-1}(X_u)\,d\mu_x(u)\,dt\,d\nu(x)\\
&=I+II.
\end{align*}
Since $X_u\subset X_x$, if $u\in X_x$, and $p\ge 1$, then 
$$
I\le\int_X\int_0^{g(x)}\mu_x^{p-1}(X_x)\mu_x(X_x)\,dt\,d\nu(x)=\int_Xg(x)\,d\nu(x).
$$ 
Since both $Sf$ and $f$ are $\prec$-decreasing, and $p\ge 1$, then $g$ is $\prec$-decreasing and $\{g>t\}$ is a $\prec$-decreasing set. Also, if $u\in \{g>t\}\cap X_x$, then $X_u\subset\{g>t\}\cap X_x$, and hence 
$$
\int_{\{g>t\}\cap X_x}\mu_x^{p-1}(X_u)\,d\mu_x(u)\le \mu_x^p(\{g>t\}\cap X_x).
$$
Therefore, using the hypothesis,

\begin{align*}
II&\le\int_0^{\infty}\int_{X\setminus\{g>t\}}\mu_x^p(\{g>t\}\cap X_x)\,d\nu(x)\,dt\\
&\le C\int_0^{\infty}\int_{\{g>t\}}\,d\nu(x)\,dt=C\int_X g(x)\,d\nu(x).
\end{align*}
So, using H\"older's inequality,
$$
\Vert Sf\Vert_{L^p(d\nu)}^p\le C\int_X(Sf)^{p-1}(x)f(x)\,d\nu(x)\le C\Vert Sf\Vert_{L^p(d\nu)}^{p-1}\Vert f\Vert_{L^p(d\nu)}.
$$

\noindent
From this a priori estimate, one obtains the general result by a standard density argument.

\medskip

If $0<p<1$, and $f\in L^p_{\prec}(d\nu)$, set $D_t=\{f>t\}$, and 

$$
g_x(t)=\int_{D_t\cap X_x}d\mu_x(u).
$$

\noindent
 Then, using the embedding $L^p_{\rm dec}(t^{p-1})\hookrightarrow L^1$ (see \cite{Sa}), we have (observe that $g_x$ is a decreasing function):

\begin{align*}
\bigg(\int_X\bigg(\int_{X_x}f(u)\,d\mu_x(u)\bigg)^p\, d\nu(x)\bigg)^{1/p}&=\bigg(\int_X\bigg(\int_0^{\infty}g_x(t)\,dt\bigg)^p\, d\nu(x)\bigg)^{1/p}\\
&\le C \bigg(\int_X\int_0^{\infty}t^{p-1}(g_x(t))^p\,dt\, d\nu(x)\bigg)^{1/p}.
\end{align*}

Since $g_x(t)\le \mu_x(X_x)=1$, then

$$
\int_0^{\infty}t^{p-1}\int_{D_t}(g_x(t))^p\, d\nu(x)\,dt\le \int_0^{\infty}t^{p-1}\int_{D_t} \, d\nu(x)\,dt=\frac1{p}\Vert f\Vert_{L^p(d\nu)}^p.
$$

On the other hand, using the hypothesis (observe that $D_t\cap X_x$ is a decreasing set), 

\begin{align*}
\int_0^{\infty}t^{p-1}\int_{X\setminus D_t}(g_x(t))^p\, d\nu(x)\,dt\le C \int_0^{\infty}t^{p-1}\int_{D_t} \, d\nu(x)\,dt=\frac C{p}\Vert f\Vert_{L^p(d\nu)}^p.
\end{align*}

Therefore, 
$$
\bigg(\int_X\bigg(\int_{X_x}f(u)\,d\mu_x(u)\bigg)^p\, d\nu(x)\bigg)^{1/p}\le C\Vert f\Vert_{L^p(d\nu)}.
$$
$\hfill\Box$\bigbreak

The following results follow easily by particularizing on each case condition (\ref{gencon}) of Theorem~\ref{genth}:
\begin{corollary}\label{cases}\ 

\begin{enumerate}
\item {\bf(Case of equality of orders $\prec$ and $\le$)}  Let $(X,\mu_x,\le)$ satisfy the conditions given above. Let $d\nu$ be a measure on $X$, and $p>0$. Then, the Hardy operator is bounded
$$S: L^p_{\le}(d\nu)\rightarrow L^p(d\nu),
$$
if and only if, there exist a constant $C>0$ such that, for all $\le$-decreasing sets $D$,
$$
\int_{X\setminus D}\mu_x^p(D\cap X_x)\,d\nu(x)\le C\nu (D).
$$
\item {\bf (Case of ${\mathbb R}_+$)}  Condition (\ref{gencon}) of Theorem~\ref{genth} is:

$$\int_r^{\infty}\big(\frac rx\big)^p\,d\nu(x)\le C\int_0^r\,d\nu(x),$$
for all $r>0$, which is $B_p$.
\item  {\bf (Case of a tree $T$)} Condition (\ref{gencon}) of Theorem~\ref{genth}  is:
$$
\sum_{x\in T\setminus D}\frac{|x\vee D|^p}{|x|^p} \nu(x)\le C\sum_{x\in D}\nu(x),
$$
where $x\vee D$ is the largest vertex in $[o,x]\cap D$.
\item {\bf (Case of ${\mathbb R}_+^2)$}  Condition (\ref{gencon}) of Theorem~\ref{genth} is:
$$\int_{{\mathbb R}_+^2\setminus D}\frac{|D_{x}|^p}{t^p}\, d\nu(x,t)\le C\int_D\,d\nu(x,t),
$$
where $D_{x}=\{t>0: (x,t)\in D\},$ and $D$ is any decreasing set (on each variable).
\end{enumerate}
\end{corollary}
\bigbreak

\begin{remark}{\rm As we have mentioned above, the case of ${\mathbb R}_+$ was first considered in \cite{AM}. $B_p$ weights are well understood and enjoy a very rich structure (see also \cite{Sa}, \cite{CaSoA}, and \cite{CRS}  for an account of $B_p$ and normability properties of weighted Lorentz spaces). 

The discrete case $\mathbb N$ is a particular case of a tree, and can be found in \cite{CRS}.
Weights for a general tree were studied, without the monotonicity condition, in \cite{EHP} and \cite{ARS}. It is easy to prove that a weight satisfying Corollary \ref{cases} (3) must necessarily be in $B_p({\mathbb N})$ (uniformly)  on each geodesic (see \cite{CRS}), but the converse is not true in general.
}\end{remark}

\medskip

We want to give a new example which shows that, even in dimension one ($X=\mathbb R_{+}$), Theorem~\ref{genth} can provide new results:

\medskip

Set $X_{n}=(n,n+1)$, $n\in\mathbb N\cup\{0\}$, and define the order $x\blacktriangleleft y$ if and only if $x,y\in X_{n}$ (i.e., $[x]=[y]$) and $x\le y$ (the usual order in $\mathbb R_+$). Define also the measures $d\mu_x(y)=dy/(x-[x])$ (in fact we only consider non-integer positive numbers). We introduce three different kinds of additional orders:

\begin{itemize}
  \item $\prec_1\ =\ \blacktriangleleft$.
  \item $x\prec_2y\Leftrightarrow [x]\le[y]$ and $x-[x]\le y-[y]$.
  \item $x\prec_3y\Leftrightarrow [x]<[y]$ or $x\blacktriangleleft y$ (i.e., $x\le y$).
\end{itemize}

It is easy to check that, in all cases, the orders $\blacktriangleleft$ and $\prec_j$, $j=1,2,3$,  satisfy the hypothesis of Theorem~\ref{genth}. Observe that, for this example (i.e., for the order $\blacktriangleleft$), the operator is of the form
$$
S_{\blacktriangleleft}f(x)=\frac 1{x-[x]}\int_{[x]}^x f(t)\,dt,
$$
which is a Hardy type operator with variable end-points (more general operators, with no monotonicity conditions on the function, have been considered in \cite{HS}). For each $j=1,2,3$ we have that the function $f$ is $\prec_j$-decreasing if:

\medskip
\noindent
{\it - Case} $\prec_1$: $f$ restricted to $(n,n+1)$ is decreasing, for every $n\in\mathbb N\cup\{0\}$. The decreasing sets are of the form $\cup_{n\ge0}(n,n+a_n)$, where $a_n\in(0,1)$.
 
 \noindent
 {\it - Case} $\prec_2$: $f$ restricted to $(n,n+1)$ is decreasing, for every $n\in\mathbb N\cup\{0\}$, and $f(x)\ge f(x+1)$, $x>0$. The decreasing sets are of the form $\cup_{n\ge0}(n,n+a_n)$, where $0< a_{n+1}< a_n<1$.

 \noindent
 {\it - Case} $\prec_3$: $f$ is decreasing for the usual order $\le$. The decreasing sets are of the form $(0,a)$, $a>0$.

 \medskip
 
 A direct application of Theorem~\ref{genth} gives:
 
 \begin{corollary}\label{threeord}
 
 Let $u$ be a weight in $\mathbb R_+$. Then,  $S_{\blacktriangleleft}:L^p_{\prec_j}(u(x)\,dx)\rightarrow L^p(u(x)\,dx)$ ($j=1,2,3$), if and only if, 
  
\medskip
\noindent
{\it Case} $\prec_1$:  For every sequence $\{a_n\}_n\subset [0,1]$, 
\begin{equation}
\label{one}
\sum_{n=0}^{\infty}\int_{a_n}^1\bigg(\frac {a_n}t\bigg)^p u(n+t)\, dt\le C\sum_{n=0}^{\infty}\int_0^{a_n} u(n+t)\,dt.
\end{equation}

\medskip
\noindent
{\it Case} $\prec_2$:   For every decreasing sequence $\{a_n\}_n\subset [0,1]$, 
\begin{equation}
\label{two}
\sum_{n=0}^{\infty}\int_{a_n}^1\bigg(\frac {a_n}t\bigg)^p u(n+t)\, dt\le C\sum_{n=0}^{\infty}\int_0^{a_n} u(n+t)\,dt.
\end{equation}

\medskip
\noindent
{\it Case} $\prec_3$:   For every $n\in\mathbb N\cup\{0\}$ and  $a\in [0,1]$, 
\begin{equation}
\label{three}
 \int_{a}^1\bigg(\frac {a}t\bigg)^p u(n+t)\, dt\le C \int_0^{n+a} u(t)\,dt.
\end{equation}
\end{corollary}

\bigbreak

\begin{remark}{\rm We observe that (\ref{one})  of Corollary~\ref{threeord} is stronger than (\ref{two}), which follows from the fact that the order $\prec_1$ is stronger than $\prec_2$. Similarly for $\prec_2$ and $\prec_3$ (in fact, to see that (\ref{three}) is weaker than (\ref{two}), given $n\in\mathbb N\cup\{0\}$ and  $a\in [0,1]$, it suffices to consider the decreasing  sequence $a_0=\cdots=a_{n-1}=1$, $a_n=a$, and $a_k=0$, $k=n+1,\cdots$).

Also, if $f$ is decreasing (for the classical order $\le$), then
$$S_{\blacktriangleleft} f(x)\le Sf(x)=\frac 1x\int_0^xf(t)\,dt.
$$
Therefore, we obtain that the boundedness of $S$ on $L^p_{\le}(u(x)\,dx)$ (i.e., the classical $B_p$ condition in \cite{AM}) is stronger than (\ref{three}), since $\prec_3\ =\ \le$. A direct argument for this fact follows from the inequality
$$
a\frac{x-[a]}x\ge a-[a],\qquad x\ge a.
$$

}\end{remark} 

\bigbreak

\section{Weights in $B_p({\mathbb R}^n_+)$}
We will now show how to apply our previous result to obtain the weighted inequalities for the multidimensional Hardy operator, acting on decreasing functions. For simplicity we will only consider the case $n=2$, the general  case being an easy extension. We first  introduce the following notations:
\begin{align*}\label{bpn}
S_1f(x,y)&=\frac1x\int_0^xf(s,y)\,ds,\qquad S_2f(x,y)=\frac1{y}\int_0^yf(x,t)\,dt,\\  Sf(x,y)&=\frac1{xy}\int_0^x\int_0^yf(s,t)\,dt\,ds=S_1(S_2f)(x,y)=S_2(S_1f)(x,y).
\end{align*}
We denote by  $D_x=\{t>0:(x,t)\in D\},$ and $D_x^y=D\cap([0,x]\times[0,y]).$ $L^p_{\rm dec}(u)$ is  the usual cone of functions in $L^p(u)$, which are decreasing on each variable. Then, we have:

\medskip

\begin{theorem}\label{iter} If $0<p<\infty$, the following are equivalent conditions:
\begin{enumerate}
\item[(a)] $S: L^p_{\rm dec}(u)\longrightarrow L^p(u)$.
\item[(b)]  There exists a constant $C>0$ such that, for every decreasing set $D$:

\begin{align}\int_{{\mathbb R}_+^2\setminus D}\frac{|D_x^y|^p}{(xy)^p}u(x,y)\,dx\,dy\le C\int_Du(x,y)\,dx\,dy.\end{align}

\item[(c)]  $S_1,S_2: L^p_{\rm dec}(u)\longrightarrow L^p(u)$.
\end{enumerate}
\end{theorem}

\medskip

\noindent{\bf Proof:} 

That (a) implies (b) follows as usual: Taking $f=\chi_D$, and using the fact that 

$$Sf(x,y)=\frac{|D_x^y|}{xy},$$
then (\ref{bpn})  is a consequence of the hypothesis.

\medskip

To show that (b) implies (c), we observe that if $(x,y)\not\in D$, then 

$$[0,x]\times D_x\subset D_x^y,
$$
and hence, $x|D_x|\le|D_x^y|$. Therefore, 

\begin{align*}
\int_{{\mathbb R}_+^2\setminus D}\frac{|D_{x}|^p}{y^p}\,u(x,y)\, dx\,dy&  \le\int_{D^c}\frac{|D_x^y|^p}{(xy)^p}u(x,y)\,dx\,dy\\
&\le C\int_Du(x,y)\,dx\,dy,
\end{align*}
and the result follows from Corollary \ref{cases} (4). Similarly for $S_1$.
\medskip

(c) implies (a):  Iterate and observe that $S_jf$ is decreasing  if $f$ is decreasing. $\hfill\Box$\bigbreak

\begin{remark}{\rm
The iteration technique used to prove Theorem~\ref{iter} can be extended very easily to other settings. For example, we could consider in $\mathbb N^2$ the operator
$$
S(\{a_{n,m}\}_{n,m})=\frac1{nm}\sum_{j=1}^n\sum_{k=1}^ma_{j,k},
$$
acting on decreasing two-indexes sequences, and obtain the characterization of the boundedness of $S$ on the weighted sequence spaces $\ell^p(\{u_{n,m}\}_{n,m})$, for general weights $\{u_{n,m}\}_{n,m}$, which improves some of the results in \cite{BHP} proved only for product weights.
}\end{remark}

Condition (\ref{bpn}) in ${\mathbb R}^n_+$ takes the following form:

$$
\int_{{\mathbb R}_+^n\setminus D}\frac{|D\cap([0,x_1]\times\cdots\times[0,x_n])|^p}{(x_1\cdots x_n)^p}u(x_1,\dots,x_n)\,dx_1\cdots dx_n\le C\int_Du(x_1,\dots,x_n)\,dx_1\cdots dx_n,
$$
which will be denoted by $u\in B_p({\mathbb R}^n_+)$. Observe that since $|D\cap([0,x_1]\times\cdots\times[0,x_n])|\le x_1\cdots x_n$, then $B_p({\mathbb R}^n_+)\subset B_q({\mathbb R}^n_+)$, if $p<q$.

\medskip

We will now prove that, as in the one-dimensional case (see \cite{AM} for the original result and \cite{Ne} for a different proof, related to the one we will use), $B_p({\mathbb R}^n_+)$ satisfies the $p-\varepsilon$ condition.

\medskip

\begin{theorem}
If $u\in B_p({\mathbb R}^n_+)$, $1\le p<\infty$, then there exists an $\varepsilon>0$ such that $u\in B_{p-\varepsilon}({\mathbb R}^n_+)$.
\end{theorem}

\medskip

\noindent{\bf Proof:}  We will only consider the case $n=2$ ($n\ge3$ follows similarly). Using Theorem \ref{iter}, it suffices to show that $S_j:L^{p-\varepsilon}_{\rm dec}(u)\rightarrow L^{p-\varepsilon}(u)$, $j=1,2$, and by symmetry, we may only consider the case $j=2$. Take any decreasing set $D\subset\mathbb R_+^2$. Then there exists a decreasing function $h:\mathbb R_+\rightarrow\mathbb R_+$ such that
$$
D=\{(s,t)\in\mathbb R_+^2; 0<t<h(s)\}.
$$
Therefore, Corollary~\ref{cases} (4) gives:
$$
\int_0^{\infty}\int_{h(s)}^{\infty}\frac {h^p(s)}{t^p}u(s,t)\,dtds\le C\int_0^{\infty}\int_0^{h(s)}u(s,t)\,dtds.
$$
With $f=\chi_D$, we have:
$$
S_2f(s,t)=\begin{cases}
     1 & \text{ if } 0<t\le h(s), \\
   \displaystyle \frac{h(s)}t  & \text{ if } 0\le h(s)<t.
\end{cases}
$$
Iterating we can prove that, for every $m\in \mathbb N$,
$$
S^m_2f(s,t)=S_2\circ{\buildrel m\over \cdots}\circ S_2f(s,t)=\begin{cases}
     1 & \text{ if } 0<t\le h(s), \\
   \displaystyle \frac{h(s)}t  \sum_{j=0}^{m-1}\frac1{j!}\log^j\frac t{h(s)}& \text{ if } 0\le h(s)<t.
\end{cases}
$$
Hence, if $h(s)<t$, we have that the following inequality follows easily
$$
(S^m_2f(s,t))^p\ge \bigg(\frac{h(s)}t\bigg)^p\frac 1{(m-1)!}\log^{m-1}\frac t{h(s)},
$$
and
$$
\int_0^{\infty}\int_{h(s)}^{\infty}\bigg(\frac{h(s)}t\bigg)^p\frac 1{(m-1)!}\log^{m-1}\frac t{h(s)} u(s,t)\,dtds\le C^m\int_0^{\infty}\int_0^{h(s)}u(s,t)\,dtds.
$$
Thus, taking $\sigma>\max(C,1/p)$, and adding up in $m$:

\begin{align*}
&\int_0^{\infty}\int_{h(s)}^{\infty}\bigg(\frac{h(s)}t\bigg)^p\sum_{m=1}^{\infty}\frac 1{\sigma^{m-1}(m-1)!}\log^{m-1}\frac t{h(s)} u(s,t)\,dtds\\
=&\int_0^{\infty}\int_{h(s)}^{\infty}\bigg(\frac{h(s)}t\bigg)^{p-1/\sigma}u(s,t)\,dtds\\
\le&\sum_{m=1}^{\infty}\bigg(\frac{C}{\sigma}\bigg)^m\int_0^{\infty}\int_0^{h(s)}u(s,t)\,dtds=C'\int_0^{\infty}\int_0^{h(s)}u(s,t)\,dtds,
\end{align*}
and the result follows with $\varepsilon=1/\sigma$.

$\hfill\Box$\bigbreak

It was proved in \cite{BaPeSo1} that for the case of the identity operator (i.e., when considering embeddings), one cannot, in general, replace the condition on all decreasing sets by  just taking rectangles of the form $[0,a_1]\times\cdots\times[0,a_n]$, $a_j>0$. However, in the case of product weights, both conditions were equivalent (see \cite[Theorem 2.5]{BaPeSo1}). We  will now show that  in this context, $u(x)=\prod_{j=1}^nu_j(x_j)\in B_p({\mathbb R}^n_+)$, factorizes very nicely as $u_j\in B_p$, for all $j\in\{1,\dots,n\}$.

\medskip

\begin{theorem}\label{decrec}  Let $u(x)=\prod_{j=1}^nu_j(x_j)\in B_p({\mathbb R}^n_+)$ be a product weight. Then the following conditions are equivalent:

\begin{enumerate}

\item[(a)]  $u\in B_p({\mathbb R}^n_+)$.

\item[(b)]  For every $a_j>0$, $j\in\{1,\dots,n\}$,

\begin{align*}
&\int_{{\mathbb R}_+^n\setminus ([0,a_1]\times\cdots\times[0,a_n])}\frac{|([0,a_1]\times\cdots\times[0,a_n])\cap([0,x_1]\times\cdots\times[0,x_n])|^p}{(x_1\cdots x_n)^p}u(x)\,dx_1\cdots dx_n\\
\le C&\int_{[0,a_1]\times\cdots\times[0,a_n]}u(x)\,dx_1\cdots dx_n.
\end{align*}

\item[(c)] $u_j\in B_p$, $j\in\{1,\dots,n\}$.
\end{enumerate}

\end{theorem}

\medskip

\noindent{\bf Proof:}  As before, and by simplicity, we will work the details only for $n=2$. 

\medskip

If $u\in B_p({\mathbb R}^2_+)$, then evaluating (\ref{bpn}) for rectangles of the form $[0,a_1]\times[0,a_2]$ we get (b).

\medskip

Assuming now that (b) holds, we evaluate the condition to obtain:

\begin{align*}
&\int_{{\mathbb R}_+^2\setminus ([0,a_1]\times[0,a_2])}\frac{|([0,a_1]\times[0,a_2])\cap([0,x_1]\times[0,x_2])|^p}{(x_1x_2)^p}u_1(x_1)u_2(x_2)\,dx_1\,dx_2\\
=\bigg(&\int_0^{a_1}u_1(x_1)\,dx_1\bigg)a_2^p\int_{a_2}^{\infty}\frac{u_{2}(x_2)}{x_2^p}\,dx_2  +
\bigg(\int_0^{a_2}u_2(x_2)\,dx_2\bigg)a_1^p\int_{a_1}^{\infty}\frac{u_{1}(x_1)}{x_1^p}\,dx_1  \\
&\qquad + (a_1a_2)^p\bigg(\int_{a_1}^{\infty}\frac{u_{1}(x_1)}{x_1^p}\,dx_1  \bigg)\bigg(\int_{a_2}^{\infty}\frac{u_{2}(x_2)}{x_2^p}\,dx_2\bigg)\\
\le C\bigg(&\int_0^{a_1}u_1(x_1)\,dx_1\bigg)\bigg(\int_0^{a_2}u_2(x_2)\,dx_2\bigg),
\end{align*}
from which we easily deduce that, for $j=1,2$,

$$
a_j^p\int_{a_j}^{\infty}\frac{u_{j}(x_j)}{x_j^p}\,dx_j\le C \int_0^{a_j}u_j(x_j)\,dx_j,
$$
and hence $u_j\in B_p$.

\medskip

Finally, iterating the one-dimensional Hardy operator, and using the fact that $u$ is a product weight, we deduce that 

$$
S: L^p_{\rm dec}(u)\longrightarrow L^p(u),
$$ 

which is (a). $\hfill\Box$\bigbreak

\begin{remark}\label{twoin}{\rm\ 
The equivalence between (c) of Theorem~\ref{decrec} and the boundedness of the Hardy operator $S: L^p_{\rm dec}(u_1u_2)\longrightarrow L^p(u_1u_2)$, for the range $p\ge 1$, was proved in \cite{BKPS}, by using an indirect argument related to the characterization of the normability property of some multidimensional analogs of the weighted Lorentz spaces (in particular this proof did not make use of the $B_p({\mathbb R}^n_+)$ condition). For the case $p=1$ one can even prove a quantitative estimate of the constant in the $B_1$ condition, namely, if we set

$$
\Vert u\Vert_{B_1(\mathbb{R}^2_+)}=\displaystyle\sup_{D\ {\rm decreasing}}\frac{\displaystyle\int_{\mathbb{R}^2_+}S\chi_D(s,t)u(s,t)\,ds\,dt}{\displaystyle\int_Du(s,t)\,ds\,dt},
$$

\noindent
then, $\Vert u_1(x_1)u_2(x_2)\Vert_{B_1(\mathbb{R}^2_+)}=\Vert u_1\Vert_{B_1}\Vert u_2\Vert_{B_1}$.

As we pointed out in Remark~\ref{twoin}, similar results to Theorem~\ref{decrec} can be obtained, for product weights, in more general settings (for example in   $\mathbb N^2$, see \cite{BHP}).

}\end{remark}
\bigbreak

\address{
\noindent
Nicola Arcozzi\\ Dept. of Mathematics\\ University of Bologna\\
40127 Bologna, ITALY \ \ \ {\sl E-mail:} 
 {\tt arcozzi@dm.unibo.it}

\medskip
\noindent
Sorina Barza\\ Dept. of Eng. Sciences, Physics and Mathematics
\\ Karlstad University\\ SE-65188 Karlstad, SWEDEN\ \ \ {\sl E-mail:} 
  {\tt sorina.barza@kau.se}
\medskip

\noindent
J.L. Garcia-Domingo\\
Economy, Mathematics \& Computers Department\\
Universitat de Vic\\E-08500 Vic, SPAIN\ \ \ {\sl E-mail:} 
 {\tt jlgarcia@uvic.es}

\medskip

\noindent
Javier Soria\\ Dept. Appl. Math. and Analysis
\\ University of Barcelona\\ E-08007 Barcelona,
 SPAIN\ \ \ {\sl E-mail:} 
 {\tt soria@ub.edu}
}

\end{document}